\documentclass{amsart}
\usepackage{amsthm}
\title{}

\newcommand{\baseRing}[1]{\ensuremath{\mathbb{#1}}}
\newcommand{\Z}{\baseRing{Z}}
\newcommand{\R}{\baseRing{R}}
\newcommand{\C}{\baseRing{C}}

\theoremstyle{plain}
\newtheorem{theorem}{Theorem}[section]
\newtheorem{lemma}[theorem]{Lemma}
\newtheorem{corollary}[theorem]{Corollary}

\newtheorem{conj}[theorem]{Conjecture}

\theoremstyle{definition}

\newtheorem{definition}[theorem]{Definition}
\newtheorem{rem}[theorem]{Remark}

\numberwithin{equation}{section}

\newcommand{\Script}[1]{\ensuremath{{\mathcal{#1}}}}

\newcommand{\UU}{\Script{U}}

\newcommand{\CC}{\Script{C}}

\renewcommand{\AA}{\Script{A}}
\newcommand{\FF}{\Script{F}}

\newcommand{\QQ}{\Script{Q}}

\newcommand{\VV}{\Script{V}}

\newcommand{\ssl}{\mathfrak{sl}}
\newcommand{\gll}{\mathfrak{gl}}

\newcommand{\ga}{\frak{a}}
\newcommand{\gkg}{\frak{g}}

\newcommand{\kk}{K\"ahler\ }

\begin{document}
\title{Mixed Lefschetz Theorems and Hodge-Riemann Bilinear Relations}
\author{Eduardo Cattani}

\address{Eduardo Cattani: Department of Mathematics
and Statistics. University
of Massachusetts. Amherst, MA 01003, USA}
\email{cattani@math.umass.edu}

\begin{abstract}

Statements analogous to the Hard Lefschetz Theorem (HLT) and the Hodge-Riemann bilinear relations (HRR) hold in a variety of contexts: they impose severe restrictions on the cohomology algebra of a smooth compact K\"ahler manifold or on the intersection cohomology of a projective toric variety; they restrict the local monodromy of a polarized variation of Hodge structure; they impose conditions on the possible $f$-vectors of convex polytopes.  While the statements of 
these theorems depend on the choice of a K\"ahler class, or its analog, there is usually a cone of possible K\"ahler classes.  It is then natural to ask whether the HLT and HRR remain true in a mixed context.  In this note we present a unified approach to proving the mixed HLT and HRR, generalizing the results obtained by \cite{cks, gromov, timorin1, timorin2, dinh-nguyen}, and proving it in new cases such as the intersection cohomology of non-rational polytopes.

\end{abstract}

\footnotetext[1]{AMS Subject Classification:
Primary 32G20, Secondary 14F43, 32Q15, 52B20. Keywords: Mixed Hard Lefschetz Theorem, Hodge-Riemann bilinear relations, Variation of Hodge Structure.}

\maketitle

\section{Introduction}

The cohomology  of a smooth compact \kk manifold $X$ is constrained by the existence of a Hodge decomposition in each degree
\begin{equation}\label{hodge}
H^d(X,\C) \ = \ \bigoplus_{p+q = d} H^{p,q}(X)\,; \quad 
\overline{H^{p,q}(X)} = H^{q,p}(X)\,,
\end{equation}
where, in de Rham terms, $H^{p,q}(X)$ may be characterized as those
cohomology classes with a representative of bidegree $(p,q)$, and by
the existence of 
a polarized Lefschetz action on the total cohomology space $H^*(X,\C)$.  The latter structure is encoded in the Hard Lefschetz Theorem (HLT) and the Hodge-Riemann bilinear relations (HRR) (see, for example, \cite{gh}):  

\begin{theorem}[{\bf HLT}]\label{hlt1}
Let $X$ be a smooth, compact, $k$-dimensional \kk manifold  and
let $\omega \in H^{1,1}(X) \cap H^2(X,\R)$ be a \kk class.  Let 
$L_\omega \in {\rm End}(H^*(X,\C))$ denote multiplication by 
$\omega$.  Then, for each $m$, such that $0\leq m \leq k$, the
map
\begin{equation}\label{hliso} 
L_\omega^m \colon H^{k-m}(X,\C) \to H^{k+m}(X,\C)
\end{equation}
is an isomorphism.
\end{theorem}

\begin{theorem}[{\bf HRR}]\label{hrr1}
Let $X$ be a smooth, compact, $k$-dimensional \kk manifold  and
let $\omega \in H^{1,1}(X) \cap H^2(X,\R)$ be a \kk class.  Define a
real bilinear form $Q$ on $H^*(X,\C)$ by
\begin{equation}\label{bilform} 
Q(\alpha,\beta) = (-1)^{\frac{(k-d)(k-d-1)}{2}}\ \int_X \alpha\cup \beta\ ,
\end{equation}
where $\deg(\alpha) = d$ and the integral is assumed to be zero if 
$\deg(\alpha\cup \beta) \not= 2k$.
Then
\begin{equation}\label{hrineq} 
i^{p-q} \ Q(\alpha,L_\omega^m \bar \alpha) \geq 0
\end{equation}
for any
$$\alpha \in H^{p,q}(X) \cap {\rm ker}(L_\omega^{m+1})\ ;\quad p+q = k-m\,.$$
Moreover, equality holds if and only if $\alpha =0$.
\end{theorem}

Similar statements hold for the action of the local monodromy on the general fiber of a local system underlying a polarized variation of Hodge structure of weight $k$ \cite{schmid,ck1}.  In the context of mirror symmetry this statement may be viewed as dual to that for the cohomology of a smooth, compact, \kk manifold \cite{cf2}.  

In another direction, the relation between algebraic geometry and 
the combinatorics of polytopes established by toric geometry, allows
 us to prove Stanley's conjecture for simple polytopes as a consequence of HLT and HRR for toric varieties and to 
  deduce the Alexandrov-Fenchel inequality for the mixed volume of polytopes, as well as other similar properties, from the Hodge index theorem, which is also a consequence of HLT and HRR \cite{stanley1, teissier1, teissier2, khov3}.  Combinatorial proofs of the generalized Stanley conjecture for arbitrary convex polytopes are then obtained through a generalization of the HLT and HRR to convex polytopes and to the intersection cohomology of the associated projective toric varieties.  An explicit construction for the simplicial case is due to Timorin \cite{timorin2}, while the general case was first obtained by Karu \cite{karu} (see also \cite{brion, bl1, bl2, bbfk1, bbfk2, braden}).  

The structures described by the Hard Lefschetz Theorem  and the Hodge-Riemann bilinear relations have been codified in different settings and with different names appropriate to the various contexts: polarized mixed Hodge structures split over $\R$ \cite{cks}, Lefschetz modules \cite{looijenga}, Frobenius modules \cite{cf2}, and polarized Hodge structures on cohomology algebras \cite{voisin}.  They involve 
the choice of a K\"ahler class in the classical algebro-geometric situation or an appropriate $\ssl_2$-action in the variation of Hodge structure or combinatorial settings; this choice takes place in an open cone defined by the Lefschetz property (\ref{hliso}) and the positivity condition  (\ref{hrineq}).
The need to consider the action of a family of $\ssl_2$'s first arose in connection with the study of the asymptotics of variations of Hodge structure \cite{cks2} and of $L^2$ and intersection cohomologies with values in a variation of polarized Hodge structure \cite{cks}.  In the latter work it was also pointed out that the Descent Lemma \cite[Lemma~1.16]{cks} and the Purity Theorem \cite[Corollary~1.13]{cks} had implications for  mixed Lefschetz actions on the cohomology of smooth compact \kk manifolds.  Subsequently, Gromov \cite{gromov} 
explicitly stated mixed Hodge-Riemann bilinear relations and proved them in special cases.  In the algebro-geometric case these mixed theorems may be stated as follows:

\begin{theorem}[{\bf Mixed HLT}]\label{mixed_hlt}
Let $X$ be a smooth, compact, $k$-dimensional \kk manifold  and
suppose that $\omega_1,\dots,\omega_m \in H^{1,1}(X) \cap H^2(X,\R)$,
$0\leq m \leq k$,
 are  \kk classes.  Let 
$L_{\omega_j} \in {\rm End}(H^*(X,\C))$ denote multiplication by 
$\omega_j$, $1\leq j \leq m$.  Then the
map
\begin{equation*}\label{mixed_hliso} 
L_{\omega_1}\cdots L_{\omega_m} \colon H^{k-m}(X,\C) \to H^{k+m}(X,\C)
\end{equation*}
is an isomorphism.
\end{theorem}

\begin{theorem}[{\bf Mixed HRR}]\label{mixed_hrr}
Let $X$ be a smooth, compact, $k$-dimensional \kk manifold.
Suppose that $m\leq k-2$ and that $\omega_1,\dots,\omega_m, \omega_{m+1} \in H^{1,1}(X) \cap H^2(X,\R)$
 are  \kk classes. As before, let $Q$ denote the 
 intersection form (\ref{bilform}). 
 Then if 
 $$\alpha \in H^{p,q}(X) \cap {\rm ker}(L_{\omega_1}\cdots L_{\omega_m} L_{\omega_{m+1}})\ ;\quad p+q = k-m\,,$$
 we have:
$$
i^{p-q} \ Q(\alpha,L_{\omega_1}\cdots L_{\omega_m}\, \bar \alpha) \geq 0
$$
and equality holds if and only if $\alpha =0$.
\end{theorem}

In \cite[\S 2.3.A]{gromov}, Gromov proved the above theorem in the case $p=q=1$ as a form  of the Alexandrov-Fenchel inequality for 
mixed volumes and stated it in the general $p=q$ case.
Timorin \cite{timorin1} proved the above mixed theorems in the linear algebraic context, i.e. the cohomology algebra of a torus and, more recently, Dinh and Nguyen \cite{dinh-nguyen} proved Theorems \ref{mixed_hlt} and \ref{mixed_hrr} for arbitrary smooth compact 
K\"ahkler manifolds.  In another direction, Timorin \cite{timorin2} indicated how the mixed HLT and HRR could be obtained in the case of simple polytopes.  

The purpose of this note is to show that when put in the context of polarized mixed Hodge structures split over the reals, here renamed 
{\em polarized Hodge-Lefschetz modules} to conform with more recent nomenclature, the mixed theorems stated above are an easy consequence of the Descent Lemma \cite[Lemma~1.16]{cks}. It should be noted, however, that the proof of the Descent Lemma makes use of the deep relationship between polarized Hodge-Lefschetz modules and variations of Hodge structure.  The advantage of this approach lies, however, in the fact that the notion of polarized Hodge-Lefschetz modules  
encompasses all cases where the HLT and HRR hold, and consequently one obtains a unified proof of the known mixed versions as well as proofs, in cases such as non-rational polytopes, where the mixed versions had not yet been proved.

This note is organized as follows: in \S 2 we define polarized
Hodge-Lefschetz modules and explain how this notion encodes the structure in the cohomology of smooth, compact K\"ahler manifold and in the combinatorial intersection cohomology of polytopes. In \S 3 we recall the basic results about polarized mixed Hodge structures and variations of Hodge structure and state the Descent Lemma which is the key result for our purposes.  Finally, in \S 4 we show how, in the context of polarized Hodge-Lefschetz modules the mixed HLT and HRR are immediate consequences of the Descent Lemma.

\section{Polarized Hodge-Lefschetz modules}

In this section we describe the abstract setting which encodes  HLT and HRR.  We have chosen as the core, a notion similar to that of 
Lefschetz modules \cite{looijenga}, although other similar objects could have been used.

\begin{definition}\label{lefschetzproperty}
Let $V = V_*$ be a $\Z$-graded finite-dimensional real vector space.  A linear map
$N \in {\rm End}_{-2}(V)$ of pure degree $-2$ is said to satisfy
the {\em Lefschetz property} relative to $V_*$ if and only if
\begin{equation*}\label{iso}
N^\ell \colon V_\ell \to V_{-\ell}
\end{equation*}
is an isomorphism for all $\ell\geq 0$. 
An abelian subspace $\ga\subset  {\rm End}_{-2}(V)$ is said
to satisfy the  Lefschetz property if some $N\in \ga$ does.
For $N$ satisfying the  Lefschetz property,
the {\em primitive subspace} $P_\ell(N)\subset V_\ell$ is the kernel of the map:
$$N^{\ell+1} \colon V_\ell \to V_{-\ell-2}.$$

\end{definition}

We shall denote by $Y \in {\rm End}_0(V)$ the semisimple transformation
acting by multiplication by $\ell$ on $V_\ell$.  It is well known
that the pair $\{Y,N\}$ may be extended to an $\ssl_2$-triple
$\{Y,N, N^+\}$, i.e. $N^+ \in {\rm End}_{2}(V)$ and  the following commutation relations hold:
$$[N_+,N] \ =\ Y\ ;\quad [Y,N_+] \ =\ 2 N_+ \ ;\quad [Y,N] \ =\ -2 N.$$
In other words, $\{Y,N, N^+\}$ define a representation of the 
Lie algebra $\ssl_2$ on $V$.
%
It follows from the basic structure theorem of $\ssl_2$-representations that the {\em Lefschetz decomposition} holds:

\begin{theorem}\label{lefschetzdec1} 
Let $V = V_*$ be a $\Z$-graded finite-dimensional real vector space and 
$N \in {\rm End}_{-2}(V)$ an endomorphism satisfying
the {\em Lefschetz property} relative to $V_*$.  Then, for every
$m\geq 0$,
\begin{equation}\label{lefeq1}
V_m \ =\ \left({\rm ker}( N^{m+1}) \cap V_m\right) \oplus NV_{m+2}.
\end{equation}
\end{theorem}

We recall that a Hodge structure of weight $d$ on a real vector space $H$ is a decomposition of its complexification $H_\C$:
\begin{equation}\label{hodgedec} H_\C \ =\ \bigoplus_{p+q =d} H^{p,q}
\end{equation}
such that $ \overline{H^{p,q}} = H^{q,p}.$   A Hodge structure of weight $d$ on $H$ is said to be {\em polarized} if there exists a real bilinear form $Q$ form of parity $(-1)^d$ such that the Hermitian form
$Q^h(.,.):= i^{-d} Q(.,\overline{.})$ makes the decomposition 
(\ref{hodgedec}) orthogonal and such that 
$(-1)^p Q^h$ is positive definite on $H^{p,d-p}$.  Given a real vector space $V$ and a non-degenerate real bilinear form $S$ on $V$ we denote
by ${\mathfrak o}(V,S)$ the Lie algebra of infinitesimal automorphisms
of $(V,S)$.

\begin{definition}\label{lefschetzhodge}
Let $V_*$ be a $\Z$-graded finite-dimensional real vector space, $k$ a positive integer, and $S$ a non-degenerate real bilinear form of parity $(-1)^k$.  Let $\ga \subset {\mathfrak o}_{-2}(V,S)$ be an abelian subspace and
$N_0 \in \ga$.  Then $(V_*,S,\ga,N_0)$ is said to be a {\em polarized Hodge-Lefschetz module} of weight $k$ if the following is satisfied:
\begin{enumerate}
\item There is a bigrading 
\begin{equation*}\label{bigrading}
V_\C \ =\ \bigoplus_{0\leq p,q\leq k} V^{p,q}\ ; \quad
V^{q,p} = \overline{V^{p,q}},
\end{equation*}
 such that
$$(V_\ell)_\C \ =\ \bigoplus_{p+q=\ell+k} V^{p,q}.$$
Hence, the bigrading restricts to   a Hodge structure of weight $k + \ell$ on $V_\ell$.

\item $T( V^{p,q}) \subset  V^{p-1,q-1}$ for all $T\in \ga$.

\item $N_0$ satisfies the Lefschetz property.

\item For $\ell\geq 0$, the induced Hodge structure on 
$P_\ell(N_0)\subset V_\ell$ is polarized by the form
$S_\ell(.,.) := S(.,N_0^\ell \,.)$.
\end{enumerate}
\end{definition}

Given a polarized Hodge-Lefschetz module $(V_*,S,\ga,N_0)$, we will
denote by 
$\CC = \CC(V_*,S,\ga,N_0)$  the largest convex cone in $\ga$ containing $N_0$ and
such that every element in $\CC$ has the Lefschetz property
relative to $V_*$.  Since it is easy to check that
every $N\in \CC$
polarizes the Hodge-Lefschetz module, in the sense that property 4 in Definition~\ref{lefschetzhodge} holds with $N$ replacing $N_0$, we will usually refer to $\CC$ as the {\em polarizing cone} of $(V_*,S,\ga,N_0)$.

\begin{rem}\label{remark1} It is shown in \cite[Proposition 1.6]{looijenga} that if 
$(V_*,S,\ga,N_0)$ is a polarized Hodge-Lefschetz module in the sense of the above definition then the Lie algebra $\gkg(\ga,V)$ generated by all  $\ssl_2$-triples
$\{Y,N, N^+\}$, where $N$ runs over all $N\in \ga$ satisfying the
Lefschetz property, is semisimple. Hence, a Hodge-Lefschetz module
is a Lefschetz module in the sense of \cite{looijenga}.
\end{rem}

\begin{rem}\label{remark2} Given a polarized Hodge-Lefschetz module $(V_*,S,\ga,N_0)$ we can construct two filtrations in $V_\C$:
$$ W_\ell := \bigoplus_{a\leq \ell} \,(V_a)_\C\ ;\quad F^p :=   \bigoplus_{a\geq p} \,V^{a,b}.$$
The filtration $W_*$ is increasing and defined over $\R$ while
the filtration $F^*$ is decreasing. $V_*$ is a grading of
$W_*[-k]$.  The pair $(W_*,F^*)$ defines
a mixed Hodge structure split over $\R$.  Moreover, it is polarized, in the sense of \cite{ck1},
by $(N,S)$, where $N$ is any element in the polarizing cone $\CC$ (cf. \cite[\S2]{cks}).  \end{rem}

\noindent
{\bf Example  (Cohomology of compact K\"ahler manifolds.)}
Let $X$ be a $k$-dimensional smooth compact \kk manifold and 
$V := H^*(X,\R)$.  We let $V_\ell := H^{k-l}(X,\R)$ and
$V^{p,q} := H^{k-q,k-p}(X)$. The Hodge decomposition (\ref{hodge}) implies that (2) in Definition~\ref{lefschetzhodge} is satisfied. Every cohomology class $\omega\in H^{1,1}(X) \cap H^2(X,\R)$ defines, by cup product, an element $L_\omega \in 
{\rm End}_{-2}(V_*)$ which is pure of bidegree $(-1,-1)$ relative to
the bigrading $V^{*,*}$. The Hard Lefschetz Theorem (cf. Theorem~\ref{hlt1}) asserts that if $\omega$ is a 
\kk class then $L_\omega$ satisfies the Lefschetz property.  On the other hand,  given $Q$ as in (\ref{bilform}), Theorem~\ref{hrr1} implies that the polarization condition (5)  in Definition~\ref{lefschetzhodge} is satisfied.  Hence, if we let $\ga$ denote 
$H^{1,1}(X) \cap H^2(X,\R)$ acting by multiplication on $V$ then for any \kk class $\omega$, $(V_*,Q,\ga,L_\omega)$, is a polarized Hodge-Lefschetz module of weight $k$.

\medskip

\noindent
{\bf Example  (Combinatorial intersection cohomology of polytopes.)}  Given a $k$-dimensional polytope $\Delta$ one may construct a
combinatorial intersection cohomology.  This is a real vector space
with an even grading 
$${\rm IH}(\Delta) = \bigoplus_{\ell=0}^k \, {\rm IH}^{2\ell}(\Delta)$$  
and a perfect intersection pairing 
$$S \colon {\rm IH}^{q}(\Delta) \times {\rm IH}^{2k - q}(\Delta) \to \R.$$
Setting $V_\ell =  {\rm IH}^{k-\ell}(\Delta)$ we have
$V_{\ell}=\{0\}$ for $\ell$ odd.  If we set $V^{p,p} := (V_{2p-k})_\C = ({\rm IH}^{2k-2p}(\Delta))_\C$ we obtain a mixed Hodge structure of Hodge-Tate type on $V$.  There is a natural action of the space $\ga$ of maps which are conewise linear on the normal fan of $\Delta$ and, for a strictly convex map $\psi$, the Lefschetz property is satisfied (see \cite[Theorem 0.1]{karu}).  Moreover, HRR is satisfied relative to the intersection form.  Hence, $(V_*,S,\ga,\psi)$ is a polarized Hodge-Lefschetz module of weight $k$ whose polarizing cone consists of strictly convex conewise linear maps.  We refer the reader to \cite{karu, 
bl1, bl2, bbfk1, bbfk2} for the details of this general case and describe, instead, Timorin's construction for the case of simple polytopes which is based on a beautiful
description of 
the cohomology algebra due to Pukhlikov and Khovanski\u{\i} \cite{khov2, khov1}.   We point out that,   
in the case of simple polytopes, the Lefschetz package was first obtained by McMullen \cite{mcmullen}. 

 Let us then assume that $\Delta$ is a simple $k$-dimensional polytope, i.e each vertex of $\Delta$ is incident to 
exactly $k$ facets. Let $r$ be the total number of facets of $\Delta$.  A polytope $P$ is said to be {\em analogous}
to $\Delta$ if $P$ and $\Delta$ have the same outward normal directions and if their facets are analogous when considered in a common hyperplane. Any two segments on a 
line are analogous.  The space $\CC(\Delta)$ of polytopes analogous
to $\Delta$ has a natural $\R$-cone structure under Minkowski sum.  It may be extended to a real vector space of virtual polytopes $\AA(\Delta)$ in the usual way.  The space $\AA(\Delta)$ comes equipped with a natural polynomial function $\nu$ of degree $k$ which restricts to the usual volume on $\CC(\Delta)$.  There is also a distinguished set of linear coordinates $x_1,\dots,x_r$ on $\AA(\Delta)$ defined in the following way: let $\xi_1,\dots,\xi_r$ be a choice of outward normal vectors
for $\Delta$, we then set for each $P \in \CC(\Delta)$,
$$x_i(P)\ :=\ \max_{m\in P} \langle \xi_i,m\rangle.$$

Let $\partial_i$ denote the operator $\partial/\partial x_i$ and
consider the graded algebra $\C[\partial] = \C[\partial_1,\dots,\partial_r]$ of partial differential operators on $\AA(\Delta)$ with 
constant coefficients.  Let $I \subset \C[\partial]$ denote
the ideal:
$$I\ :=\ \{D\in \C[\partial] : D\cdot \nu = 0\}.$$
Since $\nu$ is a homogeneous polynomial, $I$ is a homogeneous
ideal and the quotient algebra 
$$H(\Delta)\ :=\ \C[\partial]/I$$
is naturally graded.  We set: $V = \R[\partial]/(I\cap \R[\partial])$ 
with the grading
$$ V \ =\ \bigoplus_{j=-k}^k V_j\ ; \quad (V_{k-2\ell})_\C := (H(\Delta))_\ell\,,\  \ell=0,\dots,k.$$  
The set $\ga$ of linear operators acts on $V$ as linear transformations degree $-2$.  Moreover we have the following HLT \cite[Theorem~5.3.1]{timorin2}:

\begin{theorem}
Let $P \in \CC(\Delta)$ and $L_P = \sum x_i(P) \partial_i$.  Then
$L_P$ satisfies the Lefschetz property.
\end{theorem}

The vector space $V$ comes equipped with a non-degenerate pairing defined by $S([D_1],[D_2]) = D_1D_2\cdot\nu$ if the degrees of $D_1$
and $D_2$ are complementary and $S([D_1],[D_2]) =0$ otherwise. It then follows from \cite[Theorem~5.1.1]{timorin2} that HRR holds.  Hence,
if we define $(V_{k-2\ell})_\C := V^{k-\ell,k-\ell}$ we have that
$(V_*,S,\ga,L_P)$ is a Hodge-Lefschetz module for every $P$ in the
cone $\CC(\Delta)$.  We note that $\CC(\Delta)$ is the polarizing cone of this Hodge-Lefschetz module.

\section{Period mappings and variations of Hodge structure}

In this section we recall the basic definitions and main properties of polarized variations of Hodge structure (PVHS) and their period mappings.  We refer to \cite{periods1, periods2, periods3, luminy, schmid} for details.  We will show, in particular, that period mappings, more particularly nilpotent orbits  \cite{schmid,ck1,cks2}, may be viewed as the universal example of polarized Hodge-Lefschetz modules.  

Let $B$ be a connected complex manifold.  A {\it real
variation of Hodge structure} of weight $k$ (VHS) over $B$ is given by the data 
$(\VV, \nabla, \VV_{\R},\FF)$, where $\VV\to B$ is a
holomorphic vector bundle, $\nabla$ a flat connection on $\VV$,
$\VV_{\R}$ a flat real form, and $\FF$ a finite decreasing filtration
of $\VV$ by holomorphic subbundles ---the {\em Hodge filtration\/}---
satisfying

\begin{enumerate}

\item $\quad \nabla\FF^p \subset \Omega^1_B
\otimes \FF^{p-1}\quad\quad\quad\quad {\rm (Griffiths'
transversality)}$

\item $\quad\quad\  \VV = \FF^p \oplus \bar\FF^{k-p+1}
\quad\quad\quad {\rm (}\bar\FF = \hbox{conjugate of $\FF$ relative
to $\VV_{\R}$)}
$

\end{enumerate}
 As a $C^{\infty}$-bundle, $\VV$ may then be written as a
direct sum
\begin{equation}\label{decomposition}
\VV = \bigoplus_{p+q=k}\ \VV^{p,q}\ ,\quad \quad
\VV^{p,q} = \FF^p \cap \bar\FF^q\,;
\end{equation}
the integers $h^{p,q} = \dim\,\VV^{p,q}$ are the {\sl Hodge
numbers\/}.  A {\sl polarization\/} of the VHS is a flat
non-degenerate bilinear form $\QQ$ on $\VV$, whose specialization
at each fiber of $\VV$ polarizes the  Hodge structure induced
by (\ref{decomposition}) on the fiber.

Fixing a fiber $V_\C$ together
with the real structure $V_{\R}$, the polarizing form $Q$, the weight $k$
and the Hodge numbers $\{h^{p,q}\}$ and allowing the Hodge filtration $F$ to vary, we define the {\sl classifying space\/} $D := D(V,k,Q,\{h^{p,q}\})$ of
polarized Hodge structures.  Its Zariski closure $\check D$ in the
appropriate variety of flags consists of all filtrations $F$ in $V$
with $\dim\,F^p = \sum_{r\geq p}\,h^{r,k-r}$ satisfying
$$Q(F^p,F^{k-p+1}) = 0\ .$$
The complex Lie group $G_{\C}$ of all automorphisms of $(V,Q)$ acts 
transitively on $\check D$ ---therefore $\check D$ is smooth--- and
the group of real points $G_{\R}$ has $D$ as an open dense orbit. 
Let $\gkg \subset \gll(V)$ denote the Lie algebra of $G_{\C}$, 
$\gkg_{\R} \subset \gkg$ that of $G_{\R}$.  The choice of a base point 
$F\in \check D$ defines a filtration in $\gkg$
\begin{equation}\label{horiz}
F^a\gkg = \{\,T\in\gkg\ :\ T\,F^p \subset F^{p+a}\,\}\,.\end{equation}
The Lie algebra of the isotropy subgroup $U\subset G_{\C}$ at $F$ is 
$F^0\gkg$ and $F^{-1}\gkg/F^0\gkg$ is an ${\rm Ad}(U)$-invariant subspace
of $\gkg/F^0\gkg$.  The corresponding $G_{\C}$-invariant subbundle of the
holomorphic tangent bundle of $\check D$ is the {\sl horizontal
tangent bundle\/}, denoted by $T_h(\check D)$.  A polarized VHS over
a manifold $B$ determines ---via parallel translation to a typical
fiber---  a holomorphic map $\Phi\colon B \to D/\Gamma$ where $\Gamma$
is the monodromy group (Griffiths' period map).  By definition, it
has local liftings into $D$ whose differentials take values on the
horizontal tangent bundle.

In order to understand the local situation at infinity, we
suppose now that $B = (\Delta^*)^r$ is a product of punctured
disks, $\UU^r$ its universal cover, i.e. $\UU :=\{w\in \C : {\rm Im}(w) >0\}$, and  
$$ \Phi \colon \UU^r \to D$$
a horizontal map such that 
$$ \Phi(z + e_j) \ :=\ \gamma_j \Phi(z),$$
for some $\gamma_j \in G_\R$,
where $z=(z_1,\dots,z_r) \in \UU^r$, and $e_j$ represents the
$j$-th standard vector.  We assume that the transformations
$\gamma_j$ are unipotent.  If the period map arises from a polarized variation of Hodge structure defined over $\Z$, the Picard-Lefschetz transformations
$\gamma_j$ are automatically quasi-unipotent and passing, if necessary, to a finite cover of $B$ we may assume that they are
unipotent.  We will abuse notation and refer to such a map $\Phi$ as a 
local period map with unipotent monodromy.

We set
\begin{equation}\label{logmonodromy}
N_j \ :=\  \log \gamma_j \in F^{-1}\gkg \cap \gkg_\R,
\end{equation}  
and denote by $\ga$ the abelian subalgebra of $\gkg_\R$ generated
by $N_1,\dots,N_r$.  The following theorem follows from results in 
\cite{schmid, ck1, cks2, deligne}:

\begin{theorem}\label{pvhs} Let $\Phi \colon \UU^r \to D$ be
a local period mapping with unipotent monodromy and values in the classifying space 
$D(V,k,Q,\{h^{p,q}\})$.  Let $\ga$ be
the abelian Lie algebra generated by the logarithmic monodromies
$N_1,\dots,N_r$ and let $N_0 = N_1 + \cdots + N_r$.  Then 
there exists a grading of $V$ such that
$(V_*,Q,\ga,N_0)$ is a polarized Hodge-Lefschetz module of weight $k$.  
Moreover, every polarized Hodge-Lefschetz module $(V_*,Q,\ga,N_0)$ of
weight $k$ arises, in this manner, from a local period mapping with unipotent monodromy.
\end{theorem}

\begin{proof}
It follows from Schmid's Nilpotent Orbit Theorem \cite{schmid} that
if $\Phi \colon \UU^r \to D$ is
a local period mapping with unipotent monodromy then there exists
a filtration $F_{\rm lim} \in \check D$ such that the map
$$ (z_1,\dots,z_r) \in \UU^r \mapsto \exp \left(\sum_{j=1}^r z_j N_j
\right)\cdot F_{\rm lim}$$
takes values in $D$ for ${\rm Im}(z_j) >>0$.  We then know from
\cite{ck1} that every element in the positive cone $\CC$ spanned by $N_1,\dots, N_r$ defines the same weight filtration $W_*$ and as a consequence of Schmid's $SL_2$-orbit theorem it follows that $(W,F_{\rm lim} ,Q,N)$ is a polarized mixed Hodge structure for every $N \in \CC$.
Finally, it follows from \cite{deligne, cks2} that there exists a canonical splitting of this mixed Hodge structure over $\R$.  We thus
obtain a Hodge-Lefschetz module structure of weight $k$.

Conversely, suppose $(V_*,Q,\ga,N_0)$ is a polarized Hodge-Lefschetz module of weight $k$. Let $N_1,\dots,N_r$ be elements in the polarizing cone $\CC$ which are a basis of $\ga$ and such that $N_0 = N_1 + \cdots +N_r$.  Then as noted in Remark~\ref{remark2}, there exists filtrations $(W,F)$ defining a mixed Hodge structure 
of weight $k$ split over $\R$ and polarized by $(N_j,Q)$ for each $j=1,\dots, k$.
It then follows from \cite[Proposition~2.18]{ck1} that 
$$\exp  \left(\sum_{j=1}^k  z_j N_j\right) \cdot F \in D(V,k,Q,\{h^{p,q}\}),$$
where
$$h^{p,q} \ =\ \sum_{a=p} \dim V^{a,b}.$$
By (3) in Definition~\ref{lefschetzhodge}, the map 
$$ \Phi(z) \ =\ \exp  \left(\sum_{j=1}^k   z_j N_j\right) \cdot F$$
is horizontal and consequently it
defines a local period mapping with unipotent
monodromy, in fact a nilpotent orbit in the sense of \cite{schmid}.
\end{proof}

The following result is a restatement of the Descent Lemma \cite[Lemma~1.16]{cks}.

\begin{theorem}\label{descent}
Let $(V_*,Q,\ga,N_0)$ be a polarized Hodge-Lefschetz module of
weight $k$.  Let $T\in \ga$ be such that
$T + \lambda N_0$ has the Lefschetz property, relative to 
$V_*$, for all $\lambda >0$.  Let $\tilde V$ denote the image $T\cdot V$ graded by
$\tilde V_\ell = T\cdot V_{\ell + 1}$.  Set
$$\tilde Q(T\, u, T \, v)\  :=\ Q( u, T \, v)\ ;\quad u,v\in V.$$
The commutative subspace $\ga \subset {\mathfrak o}_{-2}(V,Q)$
acts on $\tilde V$ and we denote by 
$\tilde\ga \subset {\mathfrak o}_{-2}(\tilde V,\tilde Q)$
the induced space of endomorphisms.
Then $(\tilde V_*,\tilde Q,\tilde \ga,\tilde N_0)$ is a polarized
Hodge-Lefschetz module of weight $k-1$.
\end{theorem}

The proof of Theorem~\ref{descent} is the content of \S 2 of \cite{cks}. The argument, which makes extensive use of the equivalence between variations of Hodge structure and polarized Hodge-Lefschetz structures described by Theorem~\ref{pvhs}, is quite involved and we will not attempt to summarize it here.  We may illustrate the result by considering the case $T = N_0$.  Although this only applies to the unmixed case the discussion below makes clear the connection between the Descent Lemma and the mixed HLT and HRR.

Since $N_0 \cdot V^{0,q} = N_0 \cdot V^{p,0}  = 0$ we have:
$$ (\tilde V)_\C \ =\  N_0\cdot V_\C \ =\ \bigoplus_{1\leq p,q\leq k} N_0\cdot V^{p,q}.$$
Hence, if we set $\tilde V^{a,b} := N_0\cdot V^{a+1,b+1}$, we have a bigrading
$$ \tilde V_\C \ =\  \bigoplus_{0\leq a,b\leq k-1} \tilde V^{a,b}.$$
Moreover, setting $\tilde V_\ell := N_0 \cdot V_{\ell + 1}$ we have
$$ (\tilde V_\ell)_\C \ =\ \bigoplus_{p+q=\ell+ 1 + k} N_0\cdot V^{p,q}  \ =\ \bigoplus_{a + b=\ell+  (k-1) } \tilde V^{a,b}$$ 
and, since clearly 
$\overline{\tilde V^{a,b}} \ =\ \tilde V^{b,a}$, the bigrading 
$\tilde V^{*,*}$ defines a Hodge structure of weight $(k-1) + \ell$ on
$\tilde V_\ell$. We note that so far we have only used the fact that
$N_0$ is of pure bidegree $(-1,-1)$ relative to the bigrading  $V^{*,*}$.

In order to verify that $\tilde N_0$ satisfies the Lefschetz property, we note that it follows from the Lefschetz decomposition (\ref{lefeq1}) that
$$V_{-\ell+1} \ =\ N_0^{\ell -1}\cdot V_{\ell-1} \ =\ N_0^{\ell-1}\cdot P_{\ell-1}(N_0) \oplus N_0^\ell\cdot V_{\ell+1},$$
and consequently
$$\tilde V_{-\ell} \ =\ N_0\cdot  V_{-\ell+1}
 \ =\ N_0^{\ell+1}\cdot  V_{\ell+1}
  \ =\ N_0^\ell\cdot \tilde V_{\ell}.$$
Since clearly $N_0^\ell$ is injective in $\tilde V_\ell = N_0\cdot V_{\ell + 1}$ it follows that $\tilde N_0^\ell \colon \tilde V_\ell \to \tilde V_{-\ell}$ is an isomorphism.

It remains to check that $N_0$ polarizes the Hodge-Lefschetz structure on $\tilde V$ relative to $\tilde Q$.  We note, first of all, that it
is easy to check that if 
$Q$ is of parity $(-1)^k$ then $\tilde Q$ has parity $(-1)^{k-1}$.  
Finally, we verify the polarization condition 4 in Definition~\ref{lefschetzhodge}.  Let $u\in P_\ell(\tilde N_0) \cap \tilde V^{a,b}$ with $a+b = (k-1) + \ell$.  Then we can write $u = N_0 v$, where 
$$v \in P_{\ell +1}(N_0) \cap V^{a+1,b+1} \ ;\quad (a+1) + (b+1) = k + (\ell +1).$$
Hence
$$i^{a-b} \,\tilde Q(u,\tilde N_0^\ell\, \bar u) \ =\ 
i^{(a+1)-(b+1)} \,  Q(v,\tilde N_0^{\ell+1}\, \bar v) \geq 0$$ and equality holds if and only if $v=0$ or, equivalently, if $u=0$.

\begin{rem}\label{repeated}
Repeated application of Theorem~\ref{descent} allows us to replace
$\tilde V = T\cdot V$ by $\tilde V = T^m \cdot V$, $m\leq k$, graded
by $\tilde V_\ell = T\cdot V_{\ell + m}$.  Similarly we may replace
$\tilde V =T^m.V$ by $\tilde V = V/{\rm ker}(T^m)$.  
\end{rem}

\begin{corollary}\label{repeated2}
Let $(V_*,Q,\ga,N_0)$ be a polarized Hodge-Lefschetz module of
weight $k$ and $\CC$ its polarizing cone.  Let $T_1,\dots,T_m \in \CC$.  Set $\tilde V := T_1\cdots T_m \cdot V$ graded by 
$\tilde V_\ell = T_1\cdots T_m\cdot V_{\ell + m}$.  Define $\tilde Q$
by:
$$\tilde Q(T_1\cdots T_m\, u, T_1\cdots T_m \, v)\  :=\ Q( u, T_1\cdots T_m \, v)\ ;\quad u,v\in V.$$
Let $\tilde \ga$ denote $\ga$ viewed as endomorphisms of $\tilde V$.
Then $(\tilde V_*,\tilde Q,\tilde \ga,\tilde N_0)$ is a polarized
Hodge-Lefschetz module of weight $k-m$.
\end{corollary}

\begin{proof} This corollary follows from repeated application
of Theorem~\ref{descent}.
\end{proof}

In \cite{cks}, the Descent Lemma appears as a step towards the proof
of  a subtler result on the mixed $\ssl_2$ action on a polarized Hodge-Lefschetz module, namely the Purity Theorem \cite[(1.13)]{cks}. Although we do not yet have a geometric or combinatorial interpretation of this result, we include its statement for the sake of completeness.

Let $(V_*,Q,\ga,N_0)$ be a polarized Hodge-Lefschetz module of
weight $k$ and $T_1,\dots,T_m$ elements in the polarizing cone $\CC$. Consider the Koszul complex $K^*$ whose terms are defined by:
$$K^p \ := \ \bigoplus_{1\leq j_1 \leq \cdots \leq j_p\leq m}\ T_{j_1}\cdots T_{j_p}\cdot V$$ and whose differentials are given by the
maps
$$(-1)^{s-1} T_{j_s} \colon T_{j_1}\cdots \widehat{T_{j_s}}\cdots T_{j_p}\cdot V \to T_{j_1}\cdots T_{j_p}\cdot V$$ 
between the summands of $K^{p-1}$ and those of $K^{p}$. Letting $W_*(V)$ be
the natural filtration of $V$ defined by the grading $V_*$, 
we filter $K^p$  by
$$W_\ell(T_{j_1}\cdots T_{j_p}\cdot V)\ :=\ T_{j_1}\cdots T_{j_p}\cdot W_{\ell + p}(V).$$
\begin{theorem}\label{purity}
The cohomology of the filtered complex $K^*$ occurs entirely in 
weight zero or less.
\end{theorem}

\begin{rem} In the context of variations of Hodge structure, the
complex $K^*$ arises as an intersection cohomology complex and
Theorem~\ref{purity} was conjectured by Deligne as an analog of
Gabber's Purity Theorem in the $\ell$-adic case.
\end{rem}

\section{Mixed Hard Lefschetz Theorem and Hodge-Rieman bilinear relations}

In this section we show how, in the context of polarized Hodge-Lefschetz modules, the mixed HLT and HRR follow from the Descent Lemma.  We begin with the
following key lemma:

\begin{lemma}
Let $(V_*,Q,\ga,N_0)$ be a polarized Hodge-Lefschetz module of
weight $k$ and $\CC$ its polarizing cone. Let $W_*$ denote the
filtration defined by the grading $V_*$. Let $T_1,\dots,T_m \in \CC$, $m\leq k$. Then
\begin{equation}\label{key}
{\rm ker}(T_1\cdots T_m) \ \subset W_{m-1}\  = \ \bigoplus_{\ell\leq m-1} V_\ell\,.
\end{equation}
\end{lemma}

\begin{proof} We prove the statement by induction on $m$. For $m=1$ the result follows from
the assumption that  $T_1$ satisfies the Lefschetz property relative to $V_*$.  Let now
$$\tilde V \ =\ V/{\rm ker}(T_2\cdots T_m) \ \cong\ T_2\cdots T_m \cdot V\,.$$
Since $(V_*,Q,\ga,T_1)$ is also a polarized Hodge-Lefschetz module with the same polarizing cone, it follows from Corollary~\ref{repeated2} that $(\tilde V_*,\tilde Q,\tilde\ga,\tilde T_1)$ is a polarized Hodge-Lefschetz module of weight $k-m+1$.  Hence,
$${\rm ker}(\tilde T_1) \subset \tilde W_0$$ which  implies that
$$ {\rm ker}(T_1\cdots T_m) \ \subset W_{m-1} + {\rm ker}(T_2\cdots T_m) \subset W_{m-1}$$
since ${\rm ker}(T_2\cdots T_m) \subset W_{m-2} \subset W_{m-1}$
by inductive hypothesis.
\end{proof}

We can now state and prove the mixed HLT for polarized Hodge-Lefschetz modules.  When applied to the cohomology algebra of a smooth, compact K\"ahler manifold it becomes Theorem~\ref{mixed_hlt}.

\begin{theorem}
Let $(V_*,Q,\ga,N_0)$ be a polarized Hodge-Lefschetz module of
weight $k$ and $\CC$ its polarizing cone.  Let $T_1,\dots,T_m \in \CC$, $m\leq k$.  Then, the map
$$T_1\cdots T_m \colon V_m \to V_{-m}$$
is an isomorphism.
\end{theorem}

\begin{proof} By Lemma~\ref{key}, the map 
$$T_1\cdots T_m \colon V_m \to V_{-m}$$
is $1:1$.  Since $\dim V_m = \dim V_{-m}$, the result follows.
\end{proof}

The following result is the mixed Lefschetz decomposition for polarized Hodge-Lefschetz modules.  It reduces to Theorem~\ref{lefschetzdec1} when all the
transformations $T_j$ agree.

\begin{theorem}
Let $(V_*,Q,\ga,N_0)$ be a polarized Hodge-Lefschetz module of
weight $k$ and $\CC$ its polarizing cone.  Let $T_1,\dots,T_m,T_{m+1} \in \CC$, $m+2\leq k$.  Then
\begin{equation}\label{lefschetzdec}V_m\ =\  \left({\rm ker}(T_1\cdots T_{m+1}) \cap V_m\right) \oplus T_{m+1}\cdot V_{m+2}.
\end{equation}

\end{theorem}

\begin{proof}
By Lemma~\ref{key}, 
$${\rm ker}(T_1\cdots T_m T_{m+1}^2) \ \subset\  W_{m+1} \ 
=\ \bigoplus_{\ell\leq m+1} V_\ell\,. $$
Hence
$$ {\rm ker}(T_1\cdots T_m T_{m+1}) \cap T_{m+1}\cdot V_{m+2} \ =\ \{0\}\,.$$
Thus, it suffices to prove that 
$$
V_m\ \subset \  \left({\rm ker}(T_1\cdots T_{m+1}) \cap V_m\right) + T_{m+1}\cdot V_{m+2}.
$$
Let 
$\tilde V \ =\ V/{\rm ker}(T_1\cdots T_m) \ \cong\ T_1\cdots T_m \cdot V\,$.  Then, since  $(\tilde V_*,\tilde Q,\tilde\ga,\tilde T_{m+1})$ is a polarized Hodge-Lefschetz module of weight $k-m$, we have
$$\tilde V_0\  =\  \left({\rm ker}(\tilde T_{m+1}) \cap \tilde V_0\right) + \tilde T_{m+1} \cdot \tilde V_2\,.$$
Hence
$$V_m \ \subset\  \left({\rm ker}(T_1\cdots T_{m+1}) \cap V_m\right)\ +\  T_{m+1}\cdot V_{m+2}\  +\  {\rm ker}(T_1\cdots T_m).$$
Since, by Lemma~\ref{key} 
$${\rm ker}(T_1\cdots T_m) \cap V_m = \{0\},$$
the Theorem follows.
\end{proof}

The following is the mixed version of the Hodge-Riemann bilinear relations for polarized Hodge-Lefschetz modules.  In the geometric case, it is Theorem~\ref{mixed_hrr}.

\begin{theorem}\label{mixedhrr}
Let $(V_*,Q,\ga,N_0)$ be a polarized Hodge-Lefschetz module of
weight $k$ and $\CC$ its polarizing cone.  Let $T_1,\dots,T_m,T_{m+1} \in \CC$, $m+2\leq k$.  Then if 
$$v\in V^{p,q} \cap {\rm ker}(T_1\cdots T_{m+1}) \ ;\quad p+q = k+m,$$we have:
$$i^{p-q}\ Q\bigl(v, T_1\cdots T_{m} \, \bar v \bigr) \geq 0$$
with equality if and only if $v=0$.
\end{theorem}

\begin{proof}
Let $\tilde V = T_1\cdots T_m V$.  Then by 
Corollary~\ref{repeated2} applied to the polarized Hodge-Lefschetz
module $(V_*,Q,\ga,T_{m+1})$ we have that 
$(\tilde V_*,\tilde Q,\tilde \ga,\tilde T_{m+1})$ is a polarized
Hodge-Lefschetz module of weight $k-m$, where
$$\tilde Q(T_1\cdots T_m\, u, T_1\cdots T_m \, v)\  :=\ Q( u, T_1\cdots T_m \, v)\ ;\quad u,v\in V.$$

Let now $v\in V^{p,q} \cap {\rm ker}(T_1\cdots T_{m+1})$.
since $0\leq p,q \leq k$ and $p+q = k+m$ we must have
$p,q\geq m$ and the image
$$T_1\cdots T_m\, v \in \tilde V^{p-m,q-m} \cap {\rm ker}(\tilde T_{m+1}) \ ;\quad (p-m)+(q-m) = k-m.$$
Hence, 
$$i^{p-q} \tilde Q(T_1\cdots T_m\, v , T_1\cdots T_m\, \bar v ) = 
i^{p-q}   Q(  v , T_1\cdots T_m\, \bar v ) \geq 0.$$
Moreover, equality holds if and only if $T_1\cdots T_m\cdot v=0$
but, by Lemma~\ref{key}, $$V_m \cap {\rm ker}(T_1\cdots T_{m}) = \{0\}$$ and the Theorem is proved.
\end{proof}
\nocite{gromov, karu, bl1, bl2, looijenga, timorin1, timorin2,
timorin3, dinh-nguyen, khov1, khov2, khov3, brion, mcmullen, stanley1,
stanley2, cks, gh, bbfk1, bbfk2, teissier1, teissier2, cf2, voisin, ck1,  cks2, schmid}



\def\cprime{$'$} \def\cprime{$'$} \def\cprime{$'$}

\end{document}